\documentstyle[12pt]{article} \oddsidemargin=0in \parskip=.7ex
\begin{document}
\newcommand{\nc}{\newcommand}  \nc{\ov}{\over} \nc{\iy}{\infty}
\nc{\al}{\alpha} \nc{\inv}{^{-1}} \nc{\ep}{\varepsilon} \nc{\ph}{\varphi}
\nc{\x}{\xi} \nc{\be}{\beta} \nc{\sqep}{\sqrt{\ep}} \nc{\bq}{\begin{equation}}
\nc{\eq}{\end{equation}} \nc{\ch}{\raisebox{.4ex}{$\chi$}} \nc{\s}{\sigma}
\renewcommand{\sp}{\vskip2ex} \nc{\noi}{\noindent} \nc{\y}{\eta}
\nc{\rdt}{\mbox{det}_{2}\,}
\nc{\twotwo}[4]{\left(\begin{array}{cc}#1&#2\\&\\#3&#4\end{array}\right)}
\nc{\R}{\Re\,} \nc{\ga}{\gamma} \nc\qed{\hfill$\Box$}

\begin{center} {\bf Wiener-Hopf determinants with Fisher-Hartwig symbols}\end{center}

\begin{center}{{\bf Estelle L. Basor}\\
{\it Department of Mathematics\\
California Polytechnic State University\\
San Luis Obispo, CA 93407, USA}}\end{center}

\begin{center}{{\bf Harold Widom}\\
{\it Department of Mathematics\\
University of California\\
Santa Cruz, CA 95064, USA}}\end{center}

\begin{abstract}

With localization techniques one can obtain general limit theorems for Toeplitz 
determinants with Fisher-Hartwig singularities from the asymptotics for any symbol 
with one singularity of general type. There exists a family of these for which 
the determinants can be evaluated explicitly and their asymptotics determined. But 
for the Wiener-Hopf analogue, although there are likely analogous localization 
techniques, there is not a single example known of a symbol with Fisher-Hartwig 
singularity for which the determinant can be evaluated explicitly. In this paper 
we determine the asymptotics of Wiener-Hopf determinants for a symbol with one 
Fisher-Hartwig singularity of general type. We do this by showing that it is 
asymptotically equal to a Toeplitz determinant with symbol having the corresponding 
singularity.

\end{abstract}

\begin{center}{\bf I. Introduction}.\end{center}

The strong Szeg\"o limit theorem states that if the symbol
$\ph$ defined on the unit
circle has a sufficiently well-behaved logarithm then the determinant of the Toeplitz
matrix
\[T_n(\ph)=(\ph_{j-k})_{j,k=0,\cdots,n-1}\]
has the asymptotic behavior
\begin{equation} \det T_n(\ph)\sim G(\ph)^n\,E(\ph)\ \ \ {\rm as}\ n\to\iy,\label{Tlim}\end{equation}
where
\[G(\ph)=e^{(\log\ph)_0},\ \ \ E(\ph)=\exp\left(\sum_{k=1}^{\iy}k\,(\log\ph)_k\,(\log\ph)_{-k}\right).\]
Here subscripts denote Fourier coefficients.

Fisher and Hartwig \cite{FH} introduced a family of symbols with singularities and conjectured
the form of the asymptotics for these. If
\[\ph_{\al,\be}(e^{i\theta})=(2-2\cos\theta)^{(\al+\be)/2}\,e^{i(\theta-\pi)(\al-\be)/2},
\ \ \ 0<\theta<2\pi\]
(this symbol is said to have a pure Fisher-Hartwig singularity), then their symbols had the
form
\[\psi(z)=\ph(z)\prod_{j=1}^N\ph_{\al_j,\be_j}(z/z_j),\]
where $\ph$ satisfies the assumption of Szeg\"o's theorem and $z_1,\cdots,z_N$ are
distinct points on the unit circle. They conjectured that for some range of the
parameters the asymptotics had the form
\[\det\,T_n(\psi)\sim G(\ph)^n\;n^{\sum\al_j\be_j}\;E(\ph,\,\al_j,\,\be_j,\,z_j)\]
where $E(\ph,\,\al_j,\,\be_j,\,z_j)$ is a constant (whose value they did not conjecture).
Due to the work of many mathematicians
the conjecture has now been proved, and the constant $E(\ph,\,\al_j,\,\be_j,\,z_j)$
determined, in great generality.
The basic condition is that $|\R(\al_j\pm\be_j)|<1$. In early work the case
of several singularities was tackled directly and the proofs were quite
involved. But later it was discovered \cite{B,BH,BoS} that one could use localization
techniques which
made it possible to prove general results if one knew the asymptotics for the
symbols $\ph_{\al,\be}$ with pure singularity. Luckily
the Toeplitz determinants in these cases could be evaluated explicitly and their
asymptotics were then a straightforward matter \cite{BS1}. (For a
detailed history, and a proof of the general result, see \cite{BS}.
Stronger results for symbols with one singularity are contained in
\cite{E}, and a status report for a generalized conjecture can be
found in \cite{E2}.)

For the Wiener-Hopf analogue the symbol $\s(\x)$ is defined on the real line and equals
1 at $\pm\iy$. The finite Wiener-Hopf operator
$W_R(\s)$ acts on $L^2(0,R)$ and is equal to the identity plus the operator
with kernel
\[\frac{1}{2\pi}\int_{-\iy}^{\iy}(\s(\x)-1)\,e^{-i\,(x-y)\,\x}\,d\x.\]
Here there is also a ``Szeg\"o theorem'' for nicely-behaved symbols:
\begin{equation} \det W_R(\s)\sim G(\s)^R\,E(\s)\ \ \ {\rm as}\ R\to\iy.\label{WHlim}\end{equation}
Now
\[G(\s)=e^{\tau(0)},\ \ \ E(\s)=\exp\left(\int_{0}^{\iy}x\,\tau(x)\,\tau(-x)\,dx\right),\]
where $\tau$ is the Fourier transform of $\log\s$,
\[\tau(x)=\frac{1}{2\pi} \int_{-\iy}^{\iy} \log\s(\x)\,e^{-ix\x}\,d\x.\]

A Wiener-Hopf symbol with pure Fisher-Hartwig singularity could be defined by
\[\s_{\al,\be}(\x)=\left(\frac{\x-0i}{\x-i}\right)^{\al}\left(\frac{\x+0i}{\x+i}\right)^{\be}.\]
(We specify the arguments of $\x\pm0i$ and $\x\pm i$ to be zero or close to it when $\x$ is
large and positive.)
This has the behavior
\[ \s_{\al,\be}(\x)\sim |\x|^{\al+\be}\,e^{\frac{1}{2}i\pi(\al-\be)\,{\rm sgn}\x}\ \ \ {\rm as}\ \x\to0.\]
In \cite{MN} there was stated an analogue of the general 
Fisher-Hartwig conjecture for symbols which have a finite number of 
singularities of the above type, with explicit values of the 
constants. These symbols have the form
\[\s_0(\x)\,\prod_r\s_{\al_r,\be_r}(\x-\x_r),\]
where $\s_0$ is a symbol for which (\ref{WHlim}) holds.
For an attempted proof,
it seems reasonable to try to do what was done in the
Toeplitz case, that is, devise the proper localization techniques and
then try to evaluate the determinants in the case of the pure
singularity. However, while it is likely that quite general
localization techniques can be developed (and, in fact, have been in
certain cases) there has never been a single example where the finite
Wiener-Hopf determinant for a singular symbol could be evaluated
explicitly. 

In this paper we shall find the Wiener-Hopf asymptotics for the pure
symbols $\s_{\al,\be}.$ A minor complication encountered for these symbols is that 
the function $\s_{\al, \be} -1$
is not necessarily in $L^{1}$ but is in $L^1+L^{2}.$ The corresponding operators
$W_{R}(\s)-I$ will be Hilbert-Schmidt and we shall use the regularized determinants 
$\,\det_{2}W_{R}(\s).$ (For a Hilbert-Schmidt
operator $K$ the regularized determinant is defined by 
$\,\det_{2}(I+K)=\det\,(I+K)\,e^{-K}.$)

Previous results concerning the regularized determinants go back to 
\cite{Sem} where the asymptotics were determined when all
$\alpha_j=0$ or all $\beta_j = 0.$  In the case of jump
discontinuities, when all $\al_{j}+\be_{j} =0$, general results were 
also known. These were obtained by (two different
 kinds of) discretization which led to Toeplitz problems. For these 
results see \cite{BW}, \cite{BSW} and especially \cite{Arch} where the 
regularized determinant asymptotics were, in a sense, settled in the 
jump case. Also, a conjecture for the $\al = \be$ case was given in 
\cite{Mik} and proved in some special case.

We determine the asymptotics for the symbols $\s_{\al,\be}$
under the condition $|\R(\al\pm\be)|<1$. This is done by showing that after normalization
it is asymptotically equal to a regularized Toeplitz
determinant.\sp

\noi{\bf Theorem}. If $|\R(\al\pm\be)|<1$ then
\[\mbox{det}_{2}\,W_R(\s_{\al,\be})/G_{2}(\s_{\al,\be})^{R}\sim
\mbox{det}_{2}\,T_n(\ph_{\al,\be})/G_{2}(\ph_{\al,\be})^{n}  \]
when $R\sim 2n\to\iy$, where
\[G_{2}(\s) = \exp\left(\frac{1}{2\pi}\int_{-\iy}^{\iy}(\log \s(\zeta) 
-\s(\zeta) +1)\,d\zeta\right)\]
and
\[G_{2}(\ph) = \exp\big((\log\ph)_{0} -\ph_{0}+1\big).\]

\sp

This will be proved by finding exact formulas for
the Toeplitz and Wiener-Hopf determinants and regularized determinants and then
showing the above quotients are 
asymptotically equal when $R\sim2n$. These formulas are expressed in terms
of Fredholm determinants of operators acting on $L^2(0,1)$
and are obtained by using an identity
of Borodin and Okounkov \cite{BO} for Toeplitz determinants with regular symbol
and its Wiener-Hopf analogue \cite{BC}. What we do is simply stated: in both cases we
introduce a parameter to regularize the symbol, apply the identity,
and then take the limit.\sp

We remark that in the Toeplitz case, computing a regularized determinant
from an ordinary determinant is easy and it causes no problem to go back
and forth between the two in computations. The same is true for the
finite Wiener-Hopf operators when $\s_{\al, \be}-1$ is in 
$L^{1}.$ This holds exactly when $\be=\al$, and the following corollary 
will follow from the theorem.\sp

\noi{\bf Corollary}. If $|\R(\al)|<1/2$ then
\[\det \,W_R(\s_{\al,\al})\sim e^{-R\al}\,\det \,T_n(\ph_{\al,\al})\]
when $R\sim 2n\to\iy$.\sp

%%%%%%%%%%%%%%%%%%%%%%%%%%%%%%%%%%%%%%%%%%

\begin{center}{\bf II. The case of {\boldmath $\be=\al$}}\end{center}

We do this case first since the general
case has extra complications, but the main points are the same and will not have to
be repeated.

To state the Borodin-Okounkov identity, let $\ph$ be a symbol with sufficiently smooth logarithm and
let
$\ph(z)=\ph^-(z)\,\ph^+(z)$
be its Wiener-Hopf factorization, so that $\ph^+$ extends to a nonzero analytic
function inside the unit circle and $\ph^-$ outside. Let
$K_n$ be the operator on $\ell^2(\{0,\,1,\cdots\})$ with matrix entries
\[K_n(i,j)=\sum_{k=0}^{\iy}(\ph^-/\ph^+)_{n+i+k+1}\;
(\ph^+/\ph^-)_{-n-j-k-1}.\]
The identity is
\begin{equation} \det \,T_n(\ph)= G(\ph)^{n}\,E(\ph)\,\det \,(I-K_n),\label{BOT}\end{equation}
where $G(\ph)$ and $E(\ph)$ are the constants appearing in (\ref{Tlim}).
For the regularized determinant we have
\begin{equation} \rdt T_n(\ph) = G_{2}(\ph)^{n}\,E(\ph)\,\det\,
(I-K_n).\label{BOT2}\end{equation}
This follows immediately from the first identity since
\[
\rdt A = \det A \;e^{-\mbox{\scriptsize tr} (A-I)}
\]
whenever $A-I$ is trace class. Notice also that we encounter the same 
expression in computing either
\[\det T_n(\ph)/G(\ph)^{n}\]
or
\[\rdt T_n(\ph)/ G_{2}(\ph)^{n},
\]
so the crucial computation is to find $E(\ph)\,\det \,(I-K_n).$ This 
will be true when we do the Wiener-Hopf analogue as well. 

If we introduce a parameter $r<1$ the symbol $\ph_{\al,\al}(z)=(1-z)^{\al}(1-z\inv)^{\al}$
becomes $\ph_r(z)=(1-rz)^{\al}(1-rz\inv)^{\al}$, a regular symbol
for which both (\ref{BOT}) and
(\ref{BOT2}) hold. The limit of $G_{2}(\ph_{r})$ as $r\to 1$
is easily seen to be $G_{2}(\ph_{\al,\al}).$ 
The constant $E(\ph_r)$ equals $(1-r^2)^{-\al^2}$. The operator $K_n$
is the product of two Hankel operators on $\ell^2(\{0,\,1,\cdots\})$ which in this case
are the same and have $i,j$ entry
\[\frac{1}{2\pi i}\int_{|z|=1}\left(\frac{1-rz}{ 1-rz\inv}\right)^{\al}z^{n+i+j}dz
=\frac{\sin\pi\al}{ \pi}\int_0^r\left(\frac{1-rx}{ r-x}\right)^{\al}x^{n+i+j+\al}dx.\]
If the operator is denoted by $H$ then $\det\,(I-K_n)=\det\,(I-H^2)$, and so we are
interested in $\det\,(I\pm H)$. We consider at first only
$\det\,(I-H)$, the Fredholm determinant of $H$.

If $\R\al<1/2$ we can write
$H=UV$, where $U:L^2(0,r)\to\ell^2(\{0,\,1,\cdots\})$ has kernel
$U(i,x)=x^{i}$ and $V:\ell^2(\{0,\,1,\cdots\})\to L^2(0,r)$ has kernel
\[V(x,i)=\frac{\sin\pi\al}{ \pi}\left(\frac{1-rx}{ r-x}\right)^{\al}x^{n+i+\al}.\]
These are both Hilbert-Schmidt (in fact trace class) and $\det\,(I-UV)=\det\,(I-VU)$.
(See, for example, Chap. 4 of \cite{G}.)
It follows that $H$ has the same Fredholm determinant as the kernel
\[\frac{\sin\pi\al}{ \pi}\left(\frac{1-rx}{ r-x}\right)^{\al}
x^{n+\al}\frac{1}{1-xy}\]
on $L^2(0,r)$. If we set $r=(1-\ep)/(1+\ep)$ and make the substitutions
\[x\to\frac{1-x}{ 1+x},\ \ y\to\frac{1-y}{ 1+y},\]
this becomes the kernel
\begin{equation}\frac{\sin\pi\al}{ \pi}\left(\frac{x+\ep}{ x-\ep}\right)^{\al}
\left(\frac{1-x}{ 1+x}\right)^{n+\al}\frac{1}{ x+y}\label{Top}\end{equation}
on $L^2(\ep,1)$. (By this we mean that the kernels represent unitarily equivalent operators
and so have the same Fredholm determinant.) We shall determine the limit of its
Fredholm determinant as $\ep\to0$.

Denote by $A$ the operator with the above kernel
and by $A_0$ the operator with kernel
\[A_0(x,y)=\frac{\sin\pi\al}{ \pi}\frac{1}{ x+y}.\]
Let $P$ be multiplication by $\ch_{[\sqep,1]}$ and $Q$
multiplication by $\ch_{[\ep,\sqep]}$. Although our operators act on $L^2(\ep,1)$ they
can be thought of in the obvious way as acting on $L^2(0,1)$. For example the kernel
of $A_0$ can be replaced by $\ch_{[\ep,1]}(x)\,A_0(x,y)\,\ch_{[\ep,1]}(y)$.

We shall show that $P(A-A_0)$ converges in trace norm as $\ep\to0$. For this, and later use,
we give an estimate for the trace norm of certain kernels.\sp

\noi{\bf Lemma 1}. The trace norm of a kernel $f(x)g(y)/(x+y)$ on $L^2(J)$, where $J
\subset(0,\iy)$, is at
most a constant depending on $a$ times the square root of
\[\int_J|f(x)|^2\frac{dx}{ x^{1+a}}\,\cdot\,\int_J|g(x)|^2\frac{dx}{ x^{1-a}}.\]
Here $a$ belongs to $(-1,1)$ but is otherwise arbitrary.\sp

\noi{\bf Proof}. If we write the kernel as
\[\int_0^{\iy}f(x)s^{a/2}e^{-sx}\,e^{-sy}s^{-a/2}g(y)ds\]
we see that it is the (operator) product of two kernels acting
between $L^2(J)$ and
$L^2(0,\iy)$. The square of the Hilbert-Schmidt norm of the first equals
\[\int_J\int_0^{\iy}|f(x)|^2s^{a}e^{-2sx}ds\,dx,\]
which is a constant depending on $a$ times $\int_J|f(x)|^2x^{-1-a}dx$.
The second is analogous.

\sp

\noi{\bf Lemma 2}. The operator $P(A-A_0)$ converges in trace norm to the operator on
$L^2(0,1)$ with kernel
\[\frac{\sin\pi\al}{ \pi}\left[\left(\frac{1-x}{ 1+x}\right)^{n+\al}-1\right]\frac{1}{ x+y}.\]

\sp

\noi{\bf Proof}. The kernel of the operator is the sine factor times
\[\ch_{[\sqep,1]}(x)\left[\left(\frac{x+\ep}{ x-\ep}\right)^{\al}
\left(\frac{1-x}{ 1+x}\right)^{n+\al}-1\right]\frac{1}{ x+y}\]
on $L^2(\ep,1)$. If we replace the first $x$-factor by 1 then the error is of the form
\[\ch_{[\sqep,1]}(x)\,\ep_1(x)\left(\frac{1-x}{ 1+x}\right)^{n+\al}\frac{1}{ x+y}\]
where $\ep_1=O(\sqep)$. If we apply Lemma 1 with small positive $a$
we find that the operator has trace norm
$O(\ep^{1/2-\delta})$ for any $\delta>0$.
After this replacement we are left with
\[\ch_{[\sqep,1]}(x)\left[\left(\frac{1-x}{ 1+x}\right)^{n+\al}-1\right]\frac{1}{ x+y},\]
so it suffices to show that
\[\left[\left(\frac{1-x}{ 1+x}\right)^{n+\al}-1\right]\frac{1}{ x+y}\]
is trace class on $L^2(0,1)$. This also follows from Lemma 1.\sp

\noi{\bf Lemma 3}. If $\al$ is sufficiently small and $\R\al\le0$ we have as $\ep\to0$
\[\det\,(I-K_n)\sim\det\,(I-{A_0}^2)\]
\[\times\det\Big(I-(I-A_0)\inv P(A-A_0)\Big)\,\det\Big(I+(I+A_0)\inv P(A-A_0)\Big)\]
\begin{equation}\times\det\Big(I-((I-A_0)\inv Q(A-A_0)\Big)\,\det\Big(I+((I+A_0)\inv Q(A-A_0)\Big).
\label{detprod}\end{equation}

\sp

\noi{\bf Proof}. Recall that $\det\,(I-K_n)=\det\,(I-H^2)$ and that $\det\,(I\pm H)=
\det\,(I\pm A)$. The operator on $L^2(0,\,1)$ with kernel $A_0(x,\,y)$ without the sine
factor has norm 1. Hence if $\al$ is sufficiently small then
$I-A_0$ is invertible for all $\ep$ and
\[\det\,(I-A)=\det\,(I-A_0)\,\det\Big(I-(I-A_0)\inv[P(A-A_0)+Q(A-A_0)]\Big).\]
The operator in the second determinant can be written
\[\Big(I-(I-A_0)\inv Q(A-A_0)\Big)\,\Big(I-(I-A_0)\inv P(A-A_0)\Big)\]
\[-(I-A_0)\inv Q(A-A_0)(I-A_0)\inv P(A-A_0).\]
Now the operators $I-A_0$ are uniformly invertible (their inverses
have bounded norms as $\ep\to0$) if $\al$ is small enough. The same is true of
the first product above, because of the sine factor in $A-A_0$. It follows that (\ref{detprod}) will be
established if
we can show that the last term above is $o_1(1)$, i.e, its trace norm is $o(1)$. (And if
the analogous statement holds for $I+A$, which it will.)
We know that $P(A-A_0)$ converges in trace norm. It follows from the uniform invertibility
of $(I-A_0)\inv$ that it converges strongly as $\ep\to0$ to the corresponding operator on
$L^2(0,1)$, where $A_0$ does not have the $\ch_{[\ep,1]}$ factors. The same is true
of $A-A_0$ when $\R\al\le 0$ because then multiplication by 
$((x+\ep)/(x-\ep))^{\al}\,\ch_{(\sqrt\ep,\,1)}(x)$
converges strongly to $I$. Hence, since $Q$ converges strongly to 0, the same is true of
$Q(A-A_0)(I-A_0)\inv$ and this together with the trace norm convergence of $P(A-A_0)$
implies that the last term above is $o_1(1)$.\sp

We now have the ingredients necessary to derive our formula for 
$\rdt T_n(\ph_{\al,\al})/G_{2}(\ph_{\al,\al})^{n}.$\sp

\noi{\bf Lemma 4}. Assume $|\R\al|<1/2$. Let $A_0$ be the operator with kernel
\[A_0(x,y)=\frac{\sin\pi\al}{ \pi}\frac{1}{ x+y},\]
$A_1$ the operator with kernel
\[A_1(x,y)=\frac{\sin\pi\al}{ \pi}\left[\left(\frac{1-x}{ 1+x}\right)^{n+\al}-1\right]\frac{1}{ x+y},\]
and $A_2$ the operator with kernel
\[A_2(x,y)=\frac{\sin\pi\al}{ \pi}\left[\left(\frac{1+x}{ 1-x}\right)^{\al}-1\right]\frac{1}{ x+y}.\]
All act on $L^2(0,1)$. Then
\[\rdt T_n(\ph_{\al,\al})/G_{2}(\ph_{\al,\al})^{n}=4^{-\al^2}\,E(1-\sin\pi\al\;{\rm sech}\pi\x)
\,E(1+\sin\pi\al\;{\rm sech}\pi\x)\]
\[\times\det \,(I-(I-A_0)\inv A_1)\,\det\, (I+(I+A_0)\inv A_1)\]
\[\times\det \,(I-(I-A_0)\inv A_2)\,\det\, (I+(I+A_0)\inv A_2),\label{detform}\]
where the $E$ factors are those of (\ref{WHlim}).

\sp

\noi{\bf Proof}. We compute the asymptotics of the right side of (\ref{detprod}) as
$\ep\to0$ ($r\to1$). If we
make the substitutions $x\to e^{-x},\ y\to e^{-y}$ then $A_0$ becomes the operator
on $L^2(0,\log\ep\inv)$ with kernel
\[\frac{\sin\pi\al}{2\pi}\;{\rm sech}(x-y)/2.\]
This is a finite Wiener-Hopf operator with symbol $\sin\pi\al\;{\rm sech}\pi\x$
so we can use (\ref{WHlim}). We have $G(1-\sin^2\pi\al\,{\rm sech}^2\pi\x)=e^{-\al^2}$.
Since $E(\ph_r)=(1-r^2)^{-\al^2}$ and $\ep \sim (1-r)/2$ we see that
$E(\ph_r)\,\ep^{\al^2}\to4^{-\al^2}$, which gives
\[\lim_{r\to1}E(\ph_r)\,\det\,(I-{A_0}^2)=4^{-\al^2}E(1-\sin\pi\al\;{\rm sech}\pi\x)
\,E(1+\sin\pi\al\;{\rm sech}\pi\x).\]
The limit of the factor
\[\det\Big(I-(I-A_0)\inv P(A-A_0)\Big)\]
in (\ref{detprod}) is the determinant of what we get if we replace $P(A-A_0)$ by the
operator $A_1$ on $L^2(0,1)$, by the established trace norm convergence. As for the
factor
\[\det\Big(I-(I-A_0)\inv Q(A-A_0)\Big),\]
notice that
under the variable changes $x\to\ep/x,\ y\to\ep/y$ the operator $Q(A-A_0)$
on $L^2(\ep,1)$ becomes the one with kernel
\[\ch_{[\sqep,1]}(x)\,\left[\left(\frac{1+x}{ 1-x}\right)^{\al}
\left(\frac{x-\ep}{ x+\ep}\right)^{n+\al}-1\right]\frac{1}{ x+y}.\]
And, in analogy with what went before, this converges
in trace norm to the kernel $A_2(x,y)$ on $L^2(0,1)$.
Under this substitution the kernel of $A_0$ is unchanged.

All this can be done with the other factors in (\ref{detprod}). This establishes the
lemma under our assumption that $\al$ is small enough.
But the identity will hold in any connected $\al$-region containing 0 in which
both sides are analytic, so where $\ph_{\al,\al}
\in L^1$ and where $I\pm A_0$ are invertible and $A_1$ and $A_2$ trace class. This holds
for $|\R\al|<1/2$.\sp

\noi{\bf Remark}. The operators $I\pm(I\pm A_0)\inv A_1$, those involving $n$ in the
statement of the lemma, are
uniformly invertible for large $n$. It suffices to show uniform invertibility for $I\pm$
the operator with kernel the sine factor times
\[\left(\frac{1-x}{ 1+x}\right)^{n+\al}\frac{1}{ x+y}.\]
If we drop the $\al$ we make an error with operator norm $o(1)$, and then
it suffices to prove the uniform invertibility of $I\pm$  the sine factor times
\[\left(\frac{1-x}{ 1+x}\right)^{n/2}\frac{1}{ x+y}\left(\frac{1-y}{ 1+y}\right)^{n/2}.\]
Since the spectrum of $1/(x+y)$ is $[0,\pi]$ and the factors have absolute
value at most 1, the spectrum of the above also lies in $[0,\pi]$. Since it is self-adjoint
the uniform invertibility follows whenever $1/\sin\pi\al\not\in[-1,1]$, and so when
$|\R\al|<1/2$.\sp

Now let us go to Wiener-Hopf. To state the analogue of the Borodin-Okounkov identity here, 
let $\s$ be a symbol with regular logarithm.
Write $\s=\s^+\s^-$ where $\s^+$ extends to be nonzero, bounded
and analytic in the upper half-plane and $\s^-$ in the lower
and let $K_R$ be the operator on $L^2(0,\iy)$ with kernel
\[K_R(x,y)=\int_{0}^{\iy}\left(\frac{\s^{-}}{\s^{+}}
-1\right)_{R+x+z}
\left(\frac{\s^{+}}{\s^{-}}-1\right)_{-R-z-y}dz.\]
(Here for notational convenience the subscripts denote Fourier transform.)
Then for a symbol $\s$ with sufficiently regular logarithm and
satisfying $\s -1 \in L^{1}$ we have
\begin{equation} \det \,W_R(\s)= G(\s)^{R}\,E(\s)\,\det \,(I-K_R).\label{BOWH}\end{equation}
The analogous formula for the regularized determinant is the following:
for a symbol with sufficiently regular logarithm and
satisfying $\s -1\in L^{2}$,
\begin{equation} \rdt W_R(\s)= G_{2}(\s)^{R} \,E(\s)\,\det \,(I-K_R).\label{BOWH2}\end{equation}
The formula in the case that $\s -1\in L^{1}$ follows from the
results in
\cite{BC}. The derivation of the formula when $\s -1\in L^2$ is
given in the Appendix.

We introduce a small parameter $\ep$ and change $\s_{\al,\al}$
to the regular symbol
\[\s_\ep(\x)=\left(\frac{\x^2+\ep^2}{ \x^2+1}\right)^{\al}.\]
As in the Toeplitz case we can compute the finite regularized Wiener-Hopf determinant
of $\s_{\al,\al}$ for fixed $R$ by letting $\ep$ tend to zero.
To see that this is so observe that since $\s_\ep(\x)-1$
tends to $\s_{\al,\al}-1$ in $L^1$ its Fourier transforms converge in $L^2$ on any 
finite interval. The corresponding kernels
$k_{\ep}(x-y)$
and $k(x-y)$  satisfy
\[\int_{0}^{R}\int_{0}^{R}|k_{\ep}(x-y) -k(x-y)|^{2}dxdy
\leq R\int_{-\iy}^{\iy}|k_{\ep}(x) -k(x)|^{2}dx,
\]
which tends to zero.
This shows that the regularized determinants are limits of ones
with smooth symbols. This also holds for the $G_{2}(\s_{\al,\al})$ 
term. (This discussion holds also for $\s_{\al,\be}$ when $\al\neq\be$; the
regularized symbols then converge in the space $L^1+L^2$ and their Fourier
transforms still converge in $L^2$ on any finite interval.
When $\al = \be$ the above remarks
hold for the ordinary determinant as well.)

The constant $E(\s_\ep)$ which arises in (\ref{BOWH}) equals $((1+\ep)^2/4\ep)^{\al^2}$
and so $E(\s_\ep)\,\ep^{\al^2}$ has the same limit $4^{-\al^2}$ as in the Toeplitz
case.
If our finite Wiener-Hopf operator acts on $L^2(0,R)$ then the operator $K_R$ in
(\ref{BOWH}) equals the square of
the Hankel operator acting on $L^2(0,\iy)$ with kernel (in variables
$s$ and $t$)
\[\frac{1}{2\pi}\int_{-\iy}^{\iy}\left[\left(\frac{\x+\ep i}{ \x-\ep i}\right)^{\al}
\left(\frac{\x-i}{ \x+i}\right)^{\al}-1\right]\,e^{i(R+s+t)\x}\,d\x\]
\[=\frac{\sin\pi\al}{\pi}\int_{\ep}^1\left(\frac{x+\ep }{ x-\ep }\right)^{\al}
\left(\frac{1-x}{ 1+x}\right)^{\al}\,e^{-(R+s+t)x}\,dx.\]
Making a switch $UV\to VU$ as before changes this to the kernel
\[\frac{\sin\pi\al}{\pi}\left(\frac{x+\ep}{ x-\ep}\right)^{\al}
\left(\frac{1-x}{ 1+x}\right)^{\al}e^{-Rx}\frac{1}{ x+y}\]
on $L^2(\ep,1)$.

This is almost exactly the same as the kernel (\ref{Top}). The only difference is
that the expression $((1-x)/(1+x))^{n+\al}$ there is replaced here by
\[\left(\frac{1-x}{ 1+x}\right)^{\al}e^{-Rx}.\]
The same argument given above for the Toeplitz case gives the following for the
Wiener-Hopf case.\sp

\noi{\bf Lemma 5}. Assume $|\R\al|<1/2$. Let $\tilde{A}_1$ be the operator on
$L^2(0,1)$ with kernel
\[\tilde{A}_1(x,y)=\frac{\sin\pi\al}{ \pi}\left[\left(\frac{1-x}{ 1+x}\right)^\al\,
e^{-Rx}-1\right]\frac{1}{ x+y}.\]
Then with $A_0$ and $A_2$ as in Lemma~4 we have
\[\rdt W_R(\s_{\al,\al})/G_{2}(\s_{\al,\al})^{R}=4^{-\al^2}\,E(1-\sin\pi\al\;{\rm sech}\pi\x)
\,E(1+\sin\pi\al\;{\rm sech}\pi\x)\]
\[\times\det \,(I-(I-A_0)\inv \tilde{A}_1)\,\det\, (I+(I+A_0)\inv \tilde{A}_1)\]
\[\times\det\, (I-(I-A_0)\inv A_2)\,\det \,(I+(I+A_0)\inv A_2).\]

\sp

We can now establish the theorem in the case $\al=\be$. It follows from Lemmas~4 and 5 
that
\[\frac{\rdt T_n(\ph_{\al,\al})/G_{2}(\ph_{\al,\al})^{n}}
{\rdt W_R(\ph_{\al,\al})/G_{2}(\s_{\al,\al})^{n}}=\frac
{\det \,(I-(I-A_0)\inv A_1)\,\det\, (I+(I+A_0)\inv A_1)}
{\det \,(I-(I-A_0)\inv \tilde A_1)\,\det\, (I+(I+A_0)\inv \tilde A_1)}.\]
If we knew that 
$A_1-\tilde A_1
=o_1(1)$ then this together with the uniform invertibility described
in the remark following Lemma~4 would
show that the determinants involving $A_1$ and $\tilde A_1$ are asymptotically equal.
Also, because when $n$ is replaced by anything asymptotic to it the asymptotics of
$T_n(\ph_{\al,\al})$ are the same, we may replace the condition on $R$ by the stronger one
$R=2n+O(1)$. We suppose first that $R=2n$ exactly. If we apply Lemma 1 with $a=1/2$, say, we see
that we have to show that
\[\int_0^1\left|\left(\frac{1-x}{1+x}\right)^n-e^{-2nx}\right|^2\frac{dx}{ x^{3/2}}=o(1)\]
as $n\to\iy$. The part of the integral where $x$ is bounded away from 0 is clearly $o(1)$.
If $x$ is small, say $x<\delta$, we make a variable change $x\to x/2n$ and
that part of the integral becomes a constant times
\[n^{1/2}\int_0^{2\delta n}\left|e^{-x+O(x^2/n)}-e^{-x}\right|^2\frac{dx}{ x^{3/2}}.\]
The integral over $x<\sqrt n$ is $O(n^{-2})$ and the integral over $x>\sqrt n$
is $O(e^{-\sqrt n})$ since the $O(x^2/n)$ term in the exponent is $<x/2$ in the range of
integration if $\delta$ is small enough. Thus
$A_1-\tilde A_1=o_1(1)$ when $R=2n$.
The error incurred when our actual $R=2n+O(1)$ is replaced by $2n$ easily seen to be $o(1)$.
This completes the proof of the theorem in the case $\al=\be$. \sp

\noi{\bf Remark}. To prove that the Toeplitz and Wiener-Hopf determinants are
asymptotically equal it was clearly not necessary to have the individual identities given
in Lemmas~4 and 5, but only the identity for the ratio given above. 
We derived the individual identities because there was
so little extra work involved.\sp

To derive the corollary, we need only observe that
\[
\rdt T_{n}(\ph_{\al,\al}) = 
\det\, T_{n}(\ph_{\al,\al})\; e^{-\mbox{\scriptsize tr} \,(T_{n}(\ph_{\al,\al})-I)},
\]
\[\rdt W_{R}(\s_{\al,\al}) = 
\det \,W_{R}(\s_{\al,\al}) \;e^{-\mbox{\scriptsize tr} \,(W_{R}(\s_{\al,\al})-I)},
\]
\[ \mbox{tr}\, (T_{n}(\ph_{\al,\al})-I) = n\,((\ph_{\al,\al})_{0}-1),\]
\[\mbox{tr} \,(W_{R}(\s_{\al,\al})-I) = R\,\left( 
\frac{1}{2\pi}\int_{-\iy}^{\iy}(\s_{\al,\al}(\zeta) -1)\,d \zeta\right),
\]
\[G(\ph_{\al,\al}) = 1\,\,\,\,\,\,\,\,\,\,
\mbox{and}\,\,\,\,\,\,\,\,\,\,
G(\s_{\al,\al}) = e^{ -\al} .
\]
\sp

%%%%%%%%%%%%%%%%%%%%%%%%%%%%%%%%%%%%%%%%%%%%%%%

\begin{center}{\bf III. General {\boldmath$\al$} and {\boldmath $\be$}}\end{center}                
Introducing a parameter $r$ now leads to the symbol $\ph_r(z)=(1-rz)^\al(1-rz\inv)^{\be}$.
We have
$E(\ph_r)=(1-r^2)^{-\al\be}.$
Thus
we are left with the contribution from
$\det\,( I - K_{n})$ just as in the previous case.
Proceeding as we did before
leads to the Fredholm determinant of the product of two different operators on
$L^2(\ep,1)$. One has kernel
\[2^{\be-\al}\ep^{\al-\be}\,\frac{\sin\pi\be}{\pi}\,\frac{(x+\ep)^\be}{ (x-\ep)^{\al}}\,
\frac{(1-x)^{n+\al}}{(1+x)^{n+\be}}\,\frac{1}{ x+y}\]
and the kernel of the other is obtained from this by interchanging $\al$ and $\be$.
For the product the factors $2^{\be-\al}\ep^{\al-\be}$ and
$2^{\al-\be}\ep^{\be-\al}$ cancel, and so we can replace them both by 1.
The Fredholm determinant of the product equals the Fredholm determinant of
the matrix kernel
\[\twotwo{0}{{\hspace{-5ex}}\frac{\sin\pi\be}{\pi}\,\frac{(x+\ep)^{\be}}{ (x-\ep)^{\al}}\,
\frac{(1-x)^{n+\al}}{(1+x)^{n+\be}}\,\frac{1}{x+y}}
{\frac{\sin\pi\al}{\pi}\,\frac{(x+\ep)^{\al}}{(x-\ep)^{\be}}\,
\frac{(1-x)^{n+\be}}{(1+x)^{n+\al}}\,\frac{1}{x+y}}{{\hspace{-5ex}}0}.\]

We proceed as before, with $A_0$ having matrix kernel
\[\twotwo{0}{{\hspace{-1ex}}\frac{\sin\pi\be}{\pi}\frac{x^{\be-\al}}{ x+y}}
{\frac{\sin\pi\al}{\pi}\frac{x^{\al-\be}}{ x+y}}{{\hspace{-1ex}}0}.\]
We assume at first that $\al$ and $\be$ are small and purely imaginary. Then
all $x$-factors which arise will be bounded. The analogue of (\ref{detprod}) is here
\[\det\,(I-K_n)\sim\det\,(I-A_0)\]
\begin{equation}\times\det\Big(I-(I-A_0)\inv P(A-A_0)\Big)\,
\det\Big(I-(I-A_0)\inv Q(A-A_0)\Big)\label{detprod2}\end{equation}
as $r\to1$ ($\ep\to0$). To see this, observe that
$P(A-A_0)$ is $o_1(1)$ plus the operator with kernel
\[\twotwo{0}{{\hspace{-5ex}}\frac{\sin\pi\be}{\pi}
\left[\frac{(1-x)^{n+\al}}{(1+x)^{n+\be}}-1\right]\frac{x^{\be-\al}}{ x+y}}
{\frac{\sin\pi\al}{\pi}\left[\frac{(1-x)^{n+\be}}{(1+x)^{n+\al}}-1\right]
\frac{x^{\al-\be}}{ x+y}}{{\hspace{-5ex}}0}.\]
This is true following the same argument as in Lemma~2
since the function 
\[\frac{(x+\ep)^{\al}}{(x 
-\ep)^{\be}} - x^{\al -\be},\] 
and also the one with $\al$
and $\be$ interchanged, are $O(\sqrt{\ep})$ for $x \in 
(\sqrt{\ep},1).$
{}From this it follows that the analogues of Lemmas~2 and 3 hold
and also (\ref{detprod2}).

Now we proceed to find the analogue of Lemma~4. 
The kernel $Q(A-A_0)$ is $o_1(1)$ plus the operator with kernel
\[\ch_{(\ep,\sqep)}(x)\,\twotwo{0}{{\hspace{-5ex}}\frac{\sin\pi\be}{\pi}\,
\left[\frac{(x+\ep)^{\be}}{ (x-\ep)^{\al}}-x^{\be-\al}\right]\,\frac{1}{ x+y}}
{\frac{\sin\pi\al}{\pi}\,\left[\frac{(x+\ep)^{\al}}{ (x-\ep)^{\be}}-
x^{\al-\be}\right]\,\frac{1}{ x+y}}{{\hspace{-5ex}}0}.\]
If we make the variable changes $x\to\ep/x,\ y\to\ep/y$ this becomes
\[\ch_{(\sqep,1)}(x)\,\twotwo{0}{{\hspace{-5ex}}\ep^{\be-\al}\frac{\sin\pi\be}{\pi}\,
\left[\frac{(1+x)^{\be}}{ (1-x)^{\al}}-1\right]\,\frac{x^{\al-\be}}{ x+y}}
{\ep^{\al-\be}\frac{\sin\pi\al}{\pi}\,\left[\frac{(1+x^{\al}}{ (1-x)^{\be}}-1\right]\,
\frac{x^{\be-\al}}{ x+y}}{{\hspace{-5ex}}0},\]
and the kernel of $A_0$ becomes
\[\twotwo{0}{{\hspace{-2ex}}\ep^{\be-\al}\frac{\sin\pi\be}{\pi}\frac{x^{\al-\be}}{ x+y}}
{\ep^{\al-\be}\frac{\sin\pi\al}{\pi}\frac{x^{\be-\al}}{ x+y}}{{\hspace{-2ex}}0}.\]
Since determinants are unchanged if we multiply all kernels on the left by
$\left(\begin{array}{cc}I&0\\0&\ep^{\be-\al}I\end{array}\right)$ and on the right by
$\left(\begin{array}{cc}I&0\\0&\ep^{\al-\be}I\end{array}\right)$
we can remove these $\ep$ factors from both kernels. Therefore for the limit of
the second determinant in (\ref{detprod2}) the operators act
on $L^2(0,1)$ and we replace the kernel of $A_0$ by
\[\twotwo{0}{{\hspace{-1ex}}\frac{\sin\pi\be}{\pi}\frac{x^{\al-\be}}{ x+y}}
{\frac{\sin\pi\al}{\pi}\frac{x^{\be-\al}}{ x+y}}{{\hspace{-1ex}}0},\]
and that of $Q(A-A_0)$ by
\[\twotwo{0}{{\hspace{-3ex}}\frac{\sin\pi\be}{\pi}\,
\left[\frac{(1+x)^{\be}}{ (1-x)^{\al}}-1\right]\,\frac{x^{\al-\be}}{ x+y}}
{\frac{\sin\pi\al}{\pi}\,\left[\frac{(1+x^{\al}}{ (1-x)^{\be}}-1\right]\,
\frac{x^{\be-\al}}{ x+y}}{{\hspace{-3ex}}0}.\]

Finally, for the limiting operators arising from the $P(A-A_0)$ term
we multiply the matrix kernels on the left by
$\left(\begin{array}{cc}x^{(\al-\be)/2}&{\hspace{-1ex}}0\\0&{\hspace{-1ex}}x^{(\be-\al)/2}\end{array}\right)$
and on the right by
$\left(\begin{array}{cc}y^{(\be-\al)/2}&{\hspace{-1ex}}0\\0&{\hspace{-1ex}}y^{(\al-\be)/2}\end{array}\right)$.
\ Thus the kernel of $A_0$ here becomes
\begin{equation} A_0(x,y)=\twotwo{0}{{\hspace{-2ex}}\frac{\sin\pi\be}{\pi}\frac{(x/y)^{(\be-\al)/2}}{ x+y}}
{\frac{\sin\pi\al}{\pi}\frac{(x/y)^{(\al-\be)/2}}{ x+y}}{{\hspace{-2ex}}0}\label{kernel1}\end{equation}
and the kernel coming from $P(A-A_0)$ becomes
\begin{equation} A_1(x,y)=\twotwo{0}{{\hspace{-5ex}}\frac{\sin\pi\be}{\pi}
\left[\frac{(1-x)^{n+\al}}{(1+x)^{n+\be}}-1\right]\frac{(x/y)^{(\be-\al)/2}}{ x+y}}
{\frac{\sin\pi\al}{\pi}\left[\frac{(1-x)^{n+\be}}{(1+x)^{n+\al}}-1\right]
\frac{(x/y)^{(\al-\be)/2}}{ x+y}}{{\hspace{-5ex}}0}.\label{kernel2}\end{equation}
Similarly, for the limiting operators arising from the $Q(A-A_0)$ term
we multiply the matrix kernels on the left by
$\left(\begin{array}{cc}x^{(\be-\al)/2}&{\hspace{-1ex}}0\\0&{\hspace{-1ex}}x^{(\al-\be)/2}\end{array}\right)$
and on the right by
$\left(\begin{array}{cc}y^{(\al-\be)/2}&{\hspace{-1ex}}0\\0&{\hspace{-1ex}}y^{(\be-\al)/2}\end{array}\right)$.
\ The kernel of $A_0$ here becomes
\begin{equation} A_0(x,y)=\twotwo{0}{{\hspace{-2ex}}\frac{\sin\pi\be}{\pi}\frac{(x/y)^{(\al-\be)/2}}{ x+y}}
{\frac{\sin\pi\al}{\pi}\frac{(x/y)^{(\be-\al)/2}}{ x+y}}{{\hspace{-2ex}}0}\label{kernel3}\end{equation}
and the kernel coming from $Q(A-A_0)$ becomes
\begin{equation} A_2(x,y)=\twotwo{0}{{\hspace{-5ex}}\frac{\sin\pi\be}{\pi}\,
\left[\frac{(1+x)^{\be}}{ (1-x)^{\al}}-1\right]\,\frac{(x/y)^{(\al-\be)/2}}{ x+y}}
{\frac{\sin\pi\al}{\pi}\,\left[\frac{(1+x^{\al}}{ (1-x)^{\be}}-1\right]\,
\frac{(x/y)^{(\be-\al)/2}}{ x+y}}{{\hspace{-5ex}}0}.\label{kernel4}\end{equation}
None of these operations change the determinants.

For the asymptotics of $\det\,(I-A_0)$ in (\ref{detprod2}) both versions (\ref{kernel1}) and
(\ref{kernel3})
give the same determinant. If we make the variable changes $x\to e^{-x},\
y\to e^{-y}$ then (\ref{kernel1}) becomes
\begin{equation}\twotwo{0}{{\hspace{-2ex}}\frac{\sin\pi\be}{ 2\pi}\frac{e^{(\be-\al)(x-y)/2}}{
\cosh (x-y)/2}}
{\frac{\sin\pi\al}{2\pi}\frac{e^{(\al-\be)(x-y)/2}}{ \cosh (x-y)/2}}{{\hspace{-2ex}}0}\label{matrixWH}\end{equation}
on $L^2(0,\iy)$. This is the kernel of a Wiener-Hopf operator with matrix symbol.

The matrix version of (\ref{WHlim}) is well known. The
constant $G(\s)$ is defined now by using the logarithm of the
determinant of the symbol and $E(\s)$ is replaced by $\det W(\s)W(\s\inv)$.
The symbol of $I$ minus the operator with kernel (\ref{matrixWH}) is
\[\tau(\x)=\twotwo{1}{{\hspace{-2ex}}-\frac{\sin\pi\be}{\cosh\pi(\x+i(\al-\be)/2)}}
{-\frac{\sin\pi\al}{\cosh\pi(\x-i(\al-\be)2)}}{{\hspace{-2ex}}1},\]
whose determinant equals
\begin{equation} 1-\frac{\sin\pi\al\,\sin\pi\be}{ \cosh^2\pi\x-\sin^2\pi(\al-\be)/2}.\label{symdet}\end{equation}
Using the matrix version of (\ref{WHlim}) we find that the $G$ factor times the Toeplitz
$E(\ph_r)$ has limit $4^{-\al\be}$ while the $E$ factor equals $\det W(\tau)W(\tau\inv)$.

If we set $C(\al,\be)=4^{-\al\be}\,\det W(\tau)W(\tau\inv)$
we see that we have arrived at the following point.\sp

\noi{\bf Lemma 6}. If $\al$ and $\be$ are sufficiently small and purely imaginary then
\[\rdt T_n(\ph_{\al,\be})/G_{2}(\ph_{\al,\be})^{n}=C(\al,\be)\,\det\,(I-(I-A_0)\inv
A_1)\,\det\,(I-(I-A_0)\inv A_2).\]
In the first determinant on the right the kernels of $A_0$ and $A_1$ are given by
(\ref{kernel1}) and (\ref{kernel2}) and in the second determinant the kernels of $A_0$ and
$A_2$ are given by (\ref{kernel3}) and (\ref{kernel4}).

\sp

This can be extended to any connected $(\al,\be)$ region where $A_1$ and $A_2$ are trace class
and the $I-A_0$ (both versions) are invertible. Requirements just for boundedness of the
operators are $\R\al,\ \R\be<1/2$ and $|\R(\al-\be)|<1$. An application of Lemma~2 shows
that these suffice
also for the entries of $A_1$ and $A_2$ to be trace class. (Under this condition
an $a$ can be found in each case such that the integrals that arise are finite.)
As for the invertibility of $I-A_0$, a little trigonometry shows that (\ref{symdet})
is nonzero if $\cos\pi(\al+\be)/2\not\in(-\iy,-1]$,
and this holds if $|\R(\al+\be)|<1$. Since this is a connected set
and the index of the determinant is zero for
$\al$ and $\be$ small it must be zero for all these $(\al,\be)$.
Hence under
this condition $I-A_0$ is Fredholm of index zero. But we do not know that it is invertible,
so here is what we do. We know that $W(\tau)$ is invertible if $W(\tau)\,W(\tau\inv)$ is.
This operator is of the form $I$ plus trace class and its determinant is an analytic
function of $\al$ and $\be$ for $|\R(\al\pm\be)|<1$. So we assume temporarily that
$(\al,\be)$ is not in the zero set of this determinant. (Notice that the determinant is
nonzero for $\al$ and $\be$ are
small enough.) The same applies to the other version of $A_0$.
Recall that we still have the extra requirement that $\R\al,\ \R\be<1/2$.

Once we have invertibility of both versions of $I-A_0$  we can proceed as before to the
Wiener-Hopf analogue. Now the only change is that the expressions involving $1\pm x$
appearing in (\ref{kernel2}) are replaced by
$(1-x)^{\al}(1+x)^{-\be}\,e^{-Rx}$ and $(1-x)^{\be}(1+x)^{-\al}\,e^{-Rx}$
so we must add the assumption that $\R\al,\ \R\be>-1/2$.
Other than this everything is as before
until when we come to uniform invertibility. After replacing
those expressions by $e^{-2nx}$ (the error in doing this being an operator with
norm $o(1)$), the problem becomes that of uniform invertibility of $I$ minus
\[\twotwo{0}{{\hspace{-3ex}}\frac{\sin\pi\be}{\pi}e^{-2nx}\frac{(x/y)^{(\be-\al)/2}}{ x+y}}
{\frac{\sin\pi\al}{\pi}e^{-2nx}\frac{(x/y)^{(\al-\be)/2}}{ x+y}}{{\hspace{-3ex}}0}\]
on $L^2(0,1)$. We see upon making the substitutions $x\to x/2n,\ y\to y/2n$
that this is equivalent to the uniform invertibility of $I$ minus
\[\twotwo{0}{{\hspace{-3ex}}\frac{\sin\pi\be}{\pi}e^{-x}\frac{(x/y)^{(\be-\al)/2}}{ x+y}}
{\frac{\sin\pi\al}{\pi}e^{-x}\frac{(x/y)^{(\al-\be)/2}}{ x+y}}{{\hspace{-3ex}}0}\]
on $L^2(0,2n)$. This is completely equivalent to the invertibility of $I$ minus
the operator with
the same kernel on $L^2(0,\iy)$, and this in turn is equivalent to the invertibility of
$I$ minus the operator $L$ with kernel
\[L(x,y)=\twotwo{0}{{\hspace{-3ex}}\frac{\sin\pi\be}{\pi}e^{-x/2}\frac{(x/y)^{(\be-\al)/2}}{ x+y}e^{-y/2}}
{\frac{\sin\pi\al}{\pi}e^{-x/2}\frac{(x/y)^{(\al-\be)/2}}{ x+y}e^{-y/2}}{{\hspace{-3ex}}0}.\]\sp

\noi{\bf Lemma 7}. The operator $I-L$ is invertible when $|\R(\al\pm\be)|<1$ and
$(\al,\be)$ does not lie in the zero set of some analytic
function which is nonzero for sufficiently small $\al$ and $\be$.

\sp

\noi{\bf Proof}. Let $P$ denote multiplication by $\ch_{[0,1]}$ and think of the kernel
$A_0$ given by (\ref{kernel1}) as acting on $L^2(0,\iy)$. We know 
that $I-PA_{0}P$ is invertible except for $(\al,\be)$ in the zero set 
of some analytic function which is nonzero for $\al$ and $\be$ small. Now we can write
\[I-L=I-PA_0P-P(L-A_0)P-PL(I-P)-(I-P)LP\]
and several applications of Lemma~2 show that the operators $P(L-A_0)P,$
$\ PL(I-P),$ and
$(I-P)LP$ are all trace class. Hence whenever $I-PA_0P$ is invertible the invertibility
of $I-L$ is equivalent to the invertibility of 
\[I-(I-PA_0P)\inv(P(L-A_0)P+PL(I-P)+(I-P)LP),\]
which in turn is equivalent to the nonvanishing of its determinant. Since this is analytic
in $\al$ and $\be$ and nonzero for small $\al$ and $\be$ the assertion follows.

\sp

Having established the necessary invertibility and uniform invertibility the
asymptotics stated in the theorem follow, as in the case $\al=\be$. But recall that we still
have two conditions beyond the hypothesis $|\R(\al\pm\be)|<1$ of the theorem:

(i) $(\al,\be)$ does not lie in the zero set of some analytic
function $F(\al,\be)$ which is nonzero for sufficiently small $\al$ and $\be$.

(ii) $|\R\al|,\ |\R\be|<1/2$.

To remove requirement (i) we use the analyticity of the regularized determinants
and geometric means, the latter also being nonzero. Suppose $\al$ and $\be$ satisfy (ii). 
The set $S$
of $\al$ for which $F(\al,\be)$ is identically 0 in $\be$ is discrete. Assume $\al\not\in
S$. Then $\{\be:F(\al,\be)=0\}$ is discrete. Choose any $\be$ with $|\R\be|<1/2$.
There is a little circle $\Gamma$ with center $\be$ such that $F(\al,\be')\neq0$
for all $\be'\in\Gamma$. We know that
\[\frac{\rdt W_R(\s_{\al,\be'})/G_{2}(\s_{\al,\be'})^{R}}{\rdt 
T_n(\ph_{\al,\be'})/G_{2}(\ph_{\al,\be'})^{n}}\to 1\]
for all $\be'\in\Gamma$, and the denominator is also nonzero for $\be'$ inside $\Gamma$
by the known asymptotics of $\det T_n(\ph_{\al,\be'})$.
The limit holds uniformly for all $\be'\in\Gamma$. Therefore the limit holds for
$\be'=\be$ as well. So the statement of the theorem holds for all $\al$ satisfying
$|\R\al|<1/2$ except for those lying in $S$. But now we can repeat
the previous argument to show that it holds for $\al\in S$ as well.

To remove requirement (ii) we do now what we did not do earlier only because it would
have made the formulas yet more complicated. We had these requirements because of the
factors involving powers of $1-x$ in (\ref{kernel4}) and the Wiener-Hopf analogue
of (\ref{kernel2}). The exponents had to have real part greater than $-1/2$ for the
$x$-factors to belong to $L^2$. What we could have done is multiply (\ref{kernel3})
and (\ref{kernel4}) on the left by
$\left(\begin{array}{cc}(1-x)^{\al/2}&{\hspace{-1ex}}0\\0&{\hspace{-1ex}}(1-x)^{\be/2}\end{array}\right)$
and on the right by
$\left(\begin{array}{cc}(1-y)^{-\al/2}&{\hspace{-1ex}}0\\0&{\hspace{-1ex}}(1-y)^{-\be/2}\end{array}\right)$,
and multiply the Wiener-Hopf analogues of (\ref{kernel1}) and (\ref{kernel2})
on the left by
$\left(\begin{array}{cc}(1-x)^{-\al/2}&{\hspace{-1ex}}0\\0&{\hspace{-1ex}}(1-x)^{-\be/2}\end{array}\right)$
and on the right by
$\left(\begin{array}{cc}(1-y)^{\al/2}&{\hspace{-1ex}}0\\0&{\hspace{-1ex}}(1-y)^{\be/2}\end{array}\right)$.
These would not have affected the determinants but the $x$-factors (and now also
the $y$-factors) belong to $L^2$ under the weaker conditions $|\R\al|,\ |\R\be|<1$.
These hold under the hypothesis of the theorem and so are not extra conditions. The
succeeding argument holds with only minor changes with these replacements.\sp

%%%%%%%%%%%%%%%%%%%%%%%%%%%%%%%%%%%%%%%%%%%%%%%%%%%%%%%%%%%%%%%%%%

\begin{center}{\bf IV. Appendix}\end{center}

The continuous analogue of the Borodin-Okounkov
identity for generalized determinants is given by the formula
\begin{equation}
\rdt W_{R}(\s) = G_{2}(\s)^{R}\,E(\s)\, \det\,(I - K_{R}) \label{app}
\end{equation}
where
the terms $G_{2}, E$ and $K_{R}$ were defined in Section 2.
We shall show the identity is valid if  $\s -1\in L^{2},$
if the Fourier transform $k$ of $\s -1$ 
is in $L^{1}$ and satisfies 
\[\int_{-\iy}^{\iy}|x|\,|k(x)|^{2} \,dx < \iy,\]
and finally if $\s$ is nonzero and has index 
zero. As we shall see, the last two assumptions imply that 
$\s = \s^{+}\s^{-}$ where $\s^{+}$ extends to be nonzero, bounded and 
analytic in the upper half-plane and $\s^{-}$ in the lower.
All factors tend to one at $\pm \iy.$ 

That (\ref{app}) is true if $\s -1\in
L^{1}\cap L^{2}$ follows from the corresponding identity
for the ordinary determinant proved in \cite{BC} and the formula
\[
\rdt A = \det A\;e^{-\mbox{\scriptsize tr} (A -I)}.
\]
To prove (\ref{app}) under the condition $\s-1\in L^2$ we shall
approximate $\s -1$ by a sequence of functions
$\s_n-1\in L^{1}\cap L^{2}$, apply the identity to each $\s_n$ and then take a limit. 
In order to guarantee
convergence of the terms appearing in the identity, we use
a Banach algebra approach. We define the algebra ${\cal{K}}$ as the 
set of all bounded functions $\psi$ whose distributional 
Fourier transform restricted to ${\bf R}-\{0\}$ is equal to a function ${\hat \psi}$ 
satisfying   
\[\|\psi\|:=\int_{-\iy}^{\iy}|{\hat \psi}(\zeta)|d \zeta +
\left(\int_{-\iy}^{\iy}|x|\,|{\hat \psi}(\zeta)|^{2}\,d \zeta\right)^{1/2} < \iy.
\]
The Fourier transform may have a delta-function summand $c\,\delta$, and the norm
on ${\cal{K}}$ is given by 
$\|\psi\|+|c|$. This is a subalgebra of $L^{\iy}$ and it is clear that the function
$\s$ is contained in the algebra. In \cite{W} it was proved that if $\psi\in
{\cal{K}} \cap L^{2}$ then there is a sequence of functions $\psi_{n}$ with 
compact support converging to $\psi$ in both the norm of $L^{2}$ and ${\cal{K}}$. This was 
actually proved for a slightly larger Banach algebra, but the proof
is the same here.
If we apply this result to the function $\s -1$ we find a sequence of 
functions $\s_{n}-1$ for which (\ref{app}) holds. 

It remains to show that each of the
various terms has the proper limit.
Since $\s_{n}-1$ converges in $L^{2}$ to $\s-1$ the terms $\rdt W_{R}(\s_{n})$ 
and $G_{2}(\s_{n})^{R}$ converge to $\rdt W_{R}(\s)$ and 
$G_{2}(\s)^{R},$ respectively. The factor $E(\s_{n})$ is equal to $\det\, 
W(\s_{n})W(\s_{n}^{-1}) = \det\,(I - H(\s_{n})H({\tilde \s_{n}}^{-1})),$
where $H(\s)$ is the Hankel operator with kernel $\hat \s (x+y)$
and $H(\tilde \s)$ is the Hankel operator with kernel $\hat \s (-x-y)$.
Now since $\s_{n}$ converge in ${\cal{K}}$ to the invertible element
$\s$ it follows that the analogous statement is true for $\s^{-1}$
and thus the corresponding Hankel operators 
converge in the Hilbert Schmidt norm. This follows immediately from 
the definition of the norm on ${\cal{K}}.$ Finally, convergence in ${\cal{K}}$
holds for the sequences $\s_{n}^{+}$ and $\s_{n}^{-}$ and
their quotients, because the well-known projections used in their definitions
are continuous in ${\cal{K}}$. This implies the convergence of the term
$\det\,(I - K_{R}(\s_{n}))$ to the corresponding $\det\,(I - K_{R})$
and completes the proof of the identity.\sp

\begin{center}{\bf Acknowledgments}\end{center}

The authors are grateful to Albrecht B\"{o}ttcher who corrected an 
error in the corollary of the main theorem in a previous version of 
this paper and who also pointed out several previous results and 
conjectures.

The first author was supported by National Science Foundation grant DMS-0200167
and the second author by grant DMS-9732687.

%%%%%%%%%%%%%%%%%%%%%%%%%%%%%%%%%%%%%%%%%%%%%%%%%%%%%%%%%%%%%%%%%%%%%%%

\end{document}